\documentclass{amsart}
\usepackage{amsfonts}
\usepackage{amstext}
\usepackage{amsmath}
\usepackage{amssymb}
\usepackage{amsthm}
\usepackage{mathpple}
\usepackage{xypic}
\def\rat{\mathbb{Q}}
\def\integ{\mathbb{Z}}
\def\real{\mathbb{R}}
\def\dieu{\mathbb{D}}
\def\ff{\mathbb{F}}
\def\nat{\mathbb{N}}

\def\aalpha{{\pmb \alpha}}
\def\mmu{{\pmb \mu}}

\def\calg{{\mathcal G}}
\def\cale{{\mathcal E}}
\def\calh{{\mathcal H}}
\def\calo{{\mathcal O}}
\def\calp{{\mathcal P}}
\def\calq{{\mathcal Q}}
\def\idm{{\frak m}}

\DeclareMathOperator{\Height}{height}

\DeclareMathOperator{\End}{End}
\DeclareMathOperator{\gal}{Gal}
\DeclareMathOperator{\aut}{Aut}

\DeclareMathOperator{\ord}{ord}
\DeclareMathOperator{\spec}{Spec}
\DeclareMathOperator{\spf}{Spf}
\DeclareMathOperator{\gl}{GL}
\DeclareMathOperator{\Frac}{Frac}
\DeclareMathOperator{\np}{NP}
\DeclareMathOperator{\Inv}{inv}
\global\let\hom\undefined
\DeclareMathOperator{\hom}{Hom}
\def\pol{{\rm pol}}
\def\univ{{\rm univ}}
\def\perf{{\rm perf}}
\def\sep{{\rm sep}}

\def\ldot{.}
\newcommand{\abs}[1]{{\left|#1\right|}}
\newcommand{\ang}[1]{{{\langle #1 \rangle}}}

\newenvironment{alphabetize}{\begin{enumerate}

}{\end{enumerate}}

\numberwithin{equation}{section}

\theoremstyle{plain}
\newtheorem{theorem}{Theorem}[section]
\newtheorem{lemma}[theorem]{Lemma}

\newtheorem{sublemma}[theorem]{Sublemma}
\theoremstyle{remark}

\newtheorem{definition}[theorem]{Definition}

\def\tensor{\otimes}
\def\cross{\times}
\newcommand{\oneover}[1]{\frac{1}{#1}}
\def\comp{\circ}
\def\iso{\cong}
\def\inv{^{-1}}
\def\units{^\cross}
\def\ndiv{{\nmid}}

\def\inject{\hookrightarrow}
\newcommand{\st}[1]{\{#1\}}
\renewcommand{\bar}[1]{{\overline{#1}}}
\newcommand{\powser}[1]{[\![#1]\!]}
\newcommand{\laurser}[1]{(\!(#1)\!)}
\def\ra{\rightarrow}
\newcommand{\til}[1]{{\widetilde{#1}}}
\def\defeq{:=}
\begin{document}

\title{Local monodromy of $p$-divisible groups}
\date{\today}
\author{Jeffrey D. Achter}
\address{Department of Mathematics, Colorado State University, Fort Collins, CO 80523} 
\email{j.achter@colostate.edu}

\author{Peter Norman}
\address{Department of Mathematics and Statistics, University of Massachusetts, Amherst, MA 01003}
\email{norman@math.umass.edu} 

\maketitle

\begin{abstract}
A $p$-divisible group over a field $K$ admits a slope decomposition;
associated to each slope $\lambda$ is an integer $m$ and a
representation $\gal(K) \ra \gl_m(D_\lambda)$, where $D_\lambda$ is
the $\rat_p$-division algebra with Brauer invariant $[\lambda]$.  We
call $m$ the multiplicity of $\lambda$ in the $p$-divisible group.
Let $G_0$ be a $p$-divisible group over a field $k$.    Suppose that $\lambda$ is not
a slope of $G_0$, but that there exists a deformation of $G$ in which
$\lambda$ appears with multiplicity one. Assume that $\lambda\not=
(s-1)/s$ for any natural number $s>1$.  We show that there exists a  
deformation $G/R$ of $G_0/k$ such that the representation $\gal(\Frac
R) \ra \gl_1(D_\lambda)$ has large image. 
\end{abstract}

\section{Introduction}

Given a rational number $\lambda \in [0,1]$, where $\lambda=r/s$ with $\gcd(r,s) = 1$,
let
 $H_\lambda$ be the $p$-divisible group defined over $\ff_p$
whose  covariant Dieudonn\'e module is generated by a single generator
$e$ satisfying the relation $(F^{s-r}-V^r)e=0.$ 
The ring of endomorphisms of $H_{\lambda}$ which are  defined  over the algebraic 
closure of  $ \ff_p$ is an order $\calo_{H_\lambda}$
in $D_\lambda$,
the
$\rat_p$-division algebra whose Brauer invariant  is the class of
$\lambda$ in $\rat/\integ$.  By a theorem of Dieudonn\'e  and Manin a
$p$-divisible group $G$ over a field $K$ is isogenous, over the
algebraic closure $\bar K$ of $K$, to a sum $\oplus_{\lambda\in\rat}
H_{\lambda, \bar K}^{m_{\lambda}}$ \cite{manin}.  The numbers $m_\lambda$ are
uniquely determined by $G$, and we say that the slope $\lambda$
appears in $G$ with multiplicity $m_\lambda$.  

Following Gross \cite{gross}, we can
associate to $G$ the $\gal(K)$ module $V^{\lambda}(G) := 
\hom(H_{\lambda, \bar K}, G_{\bar K})\otimes \rat$. This yields (see Definition \ref{defmono}) a representation
\begin{equation}
\label{diagdefmono}
\xymatrix{
\rho^{\lambda}= \rho :\gal(K) \ar[r] & \gl_{m_{\lambda}}(D_{\lambda}),
}
\end{equation}
which we call the $\lambda$-monodromy of $G$.

We say (Definition \ref{deflarge}) that the $\lambda$-monodromy is large if there is a subgroup of 
$\gl_{m_{\lambda}}(D_{\lambda})$ that has  finite index in both
$\gl_{m_{\lambda}}( \calo_{H_\lambda})$ and $\rho (\gal(K)).$

Let $k$ be an algebraically closed field of characteristic $p$ and let
$R$ be an equicharacteristic complete local domain with residue field $k$ and fraction
field $K$. Let $G_0$ be a $p$-divisible group over $k$ and $G$ a lifting
to $R$. In these circumstances  
 the representation $\rho$ has been studied extensively.
The first major work  was done by Igusa \cite{igusa}. Assume that $G_0$ 
is the $p$-divisible group of a 
supersingular elliptic curve and that $G_K$ has slopes $0$ and $1$. 
Then the action of $\gal(K)$ on
$\hom(H_0, G_K)\iso\integ_p$ is cofinite 
in $\aut(\integ_p)=\integ_p\units$. 
Both Gross \cite{gross} and Chai \cite{chaimono} studied the case when
 $G_0$ has a single slope  $c/(c+1)$ and  $G_K$ has slope $g/(g+1)$
 with multiplicity one and  slope $1$ with multiplicity $c-g$.
 Gross showed that the  image of $\gal(K)$ in 
$\aut(H_{(g-1)/g})$ is all of $\aut(H_{(g-1)/g})$, while 
Chai showed that the slope-1 representation acting on $\rat_p^{(c-g)}$
is irreducible.  We make no attempt at a comprehensive history of this
problem, but do note that a  global analogue of this question for 
the generic Newton stratum of certain moduli spaces of PEL type has
been resolved by
Deligne and Ribet \cite{deligneribet} and Hida \cite{hidabourb}.

\begin{definition}
\label{defattainable}
Let $G_0$ be a $p$-divisible group over an algebraically closed field
$k$.  We say a slope $\lambda$ is {\em attainable} from $G_0$ if there
exists a deformation $G/R$ of $G_0$ to a complete local domain $R$
with fraction field $K$ such that $\lambda$ appears as a slope of
$G_K$ with multiplicity one; and, furthermore, we require that if
$\lambda'$ is a slope of $G_K$ with $\lambda' < \lambda$, then
$\lambda'$ appears in $G_0$ with the same multiplicity as in $G_K$.
We say that such a $G$ attains $\lambda$.
\end{definition}
(The $\lambda$ which are attainable from $G_0$ are determined
completely by the Newton polygon of $G_0$, see Theorem \ref{grothconj}.)

Our main result is:
 \begin{theorem}
\label{thmainintro}
 Let $G_0$ be a $p$-divisible group over
an algebraically closed field of characteristic $p$. Assume that a rational 
number $\lambda \in [0, 1]$ is not a slope of $G_0$, that $\lambda\not
= (s-1)/s$ for any natural number $s\ge 2$, and that
$\lambda$ is attainable from $G_0$. Then there exists a deformation
of $G_0$ which attains $\lambda$ and has   large $\lambda$-monodromy.
\end{theorem}

For applications to families of abelian varieties, we need a variant
of \ref{thmainintro} adapted to deformations of 
$p$-divisible groups equipped with
quasi-polarizations.
A principal  quasi-polarization of a $p$-divisible group 
is a self-dual isomorphism
$\Phi: G \rightarrow G^t$. 

\begin{definition}
\label{defsymattainable}
Given a principally quasi-polarized  
(or pqp for short) $p$-divisible group $(G_0/k, \Phi_0)$, $k$
algebraically closed, we say 
a rational number $\lambda \in [0,1]$ is symmetrically attainable
from $(G_0, \Phi_0)$ if there is a pqp deformation 
of $(G_0, \Phi_0)$ to a pqp $p$-divisible group
$(G, \Phi)$ over a complete local domain $R$ such that $G$ attains $\lambda$.
In this case we say that 
$(G,\Phi)$ symmetrically attains $\lambda$.
\end{definition}

\begin{theorem}
\label{thmainintropqp}
Let $(G_0, \Phi_0)$ be a pqp $p$-divisible group over an algebraically
closed field of characteristic $p$.  Assume that
$\lambda$ is not a slope of $G_0$, that $\lambda \ne (s-1)/s $ 
for any natural number $s \ge 2$,
and that $\lambda$ is
 symmetrically attainable from $(G_0, \Phi_0)$. 
Then there exists 
a pqp deformation of $(G_0, \Phi_0)$ to $(G, \Phi)$ over a 
complete local domain R
with fraction field $K$ so that $G_K$ 
symmetrically attains $\lambda$  and has large $\lambda$-monodromy.
\end{theorem}

We begin by reviewing some  facts in section two  about $p$-divisible groups, 
their automorphisms, Newton polygons, and monodromy. 
The heart of the paper is section three. There we consider
a local $p$-divisible group $G_0$ over an algebraically closed field $k$
 with $a$-number equal to one, i.e., $\dim(\hom_k(\aalpha_p, G_0))=1$. 
We also assume that $\lambda$ is positive and 
 is strictly less than any slope of $G_0$.
We construct,  using techniques of Oort, and then 
analyze a particular deformation of $G_0$
that attains $\lambda$ with multiplicity one and show it has 
large $\lambda$-monodromy. In section four we prove two 
reduction steps that allow us to complete the proof of Theorem \ref{thmainintro} 
for positive slopes, and separately we prove the case of slope zero.
The last section is devoted to proving Theorem \ref{thmainintropqp},  the 
analogue of Theorem  \ref{thmainintro} for
principally quasi-polarized $p$-divisible groups.

We thank the referee for helpful comments, and C.-L. Chai for his
suggestion to use monodromy to study $p$-divisible groups.

\section{Background on $p$-divisible groups}

\subsection{Slopes, the slope filtration and Newton polygons}
We describe the  Newton polygon of a  $p$-divisible group over a field. 
The Newton polygon determines and is determined by the isogeny class of a
$p$-divisible group over an algebraically closed field.

Let $G$ be a $p$-divisible group over a field $K$ of characteristic $p$.  By a theorem of
Grothendieck (proved in \cite{zinksf}), $G$ has a filtration 
\begin{equation}
\label{eqslopefilt}
0 = G_0 \subset G_1 \subset G_2 \subset \cdots \subset G_n = G
\end{equation}
by $p$-divisible groups so that $G_i/G_{i-1}$ is a $p$-divisible group
isogenous over the algebraic closure of $K$ to 
the direct sum  of $m_{\lambda_i}$ copies of $H_{\lambda_i}$,
 and the rational numbers $\lambda_i$
satisfy $\lambda_i < \lambda_{i+1}$.  Write $\lambda_i = r_i/s_i$
where $\gcd(r_i,s_i) = 1$.  
Then the height of successive steps in the filtration is  $\Height(G_i/G_{i-1}) = m_{\lambda_i}s_i$.
By slope we mean the slope of the Frobenius acting on the covariant 
Dieudonn\'e module, which corresponds to the Verschiebung operator
of a $p$-divisible group. For example, $\mmu_{p^{\infty}}$ has slope zero. 

The Newton polygon $\np(G)$ of $G/K$ is the convex hull of the set
\[
\st{(0,0)} \cup_{i=1}^n \st{(\Height(G_i), \sum_{j=1}^i
\Height(G_i/G_{i-1}) \lambda_i)}\subset \integ^2 \subset \real^2.
\]
It is a lower-convex polygon connecting $(0,0)$ to $(\Height(G),
\dim(G))$ with slopes in the interval $[0,1]$ and integral 
breakpoints.  Any such
Newton polygon is actually realized as the Newton polygon of a
$p$-divisible group.

Let $\nu_1$ and $\nu_2$ be two Newton polygons.  Following Oort 
we say that
$\nu_1 \succeq \nu_2$ if $\nu_1$ and $\nu_2$ share the same
endpoints and if every point of $\nu_1$ is on or below that of
$\nu_2$.

Let $G_0/k$ be a $p$-divisible group over an algebraically closed  field $k$.
 By a deformation of $G_0$ we mean a
$p$-divisible group $G$ over a local ring $R$ equipped with
isomorphisms $R/\idm_R \iso
k$ and $G\cross R/\idm_R \iso G_0/k$.

Grothendieck \cite[Theorem 2.1.3]{katzsf} proved that the Newton 
polygon goes up  under
specialization and conjectured, conversely, that one can always
achieve an arbitrary ``lower'' Newton polygon by deformation.  Oort
\cite{oortnpfg} proved this converse, and we make use of his ideas
throughout this paper.

\begin{theorem}\cite{oortnpfg}\label{grothconj} Let $G_0$ be a $p$-divisible group over 
an algebraically closed field $k$.  Let 
$\nu$ be an arbitrary Newton polygon.  Then there exists a
deformation $G/R$ of $G_0$ such that $\np(G\cross \Frac R) = \nu$ if and
only if $\np(G_0) \preceq \nu$.
\end{theorem}

There is a variant of this theory for principally quasi-polarized
$p$-divisible groups.  A Newton polygon is called symmetric if, in the
notation introduced above, $\lambda_i = 1 - \lambda_{n+1-i}$ and
$m_{\lambda_i} = m_{\lambda_{n+1-i}}$.
  Much in the vein of Theorem \ref{grothconj} Oort proves
that if $G_0$ is a principally quasi-polarized $p$-divisible group,
and 
if $\nu$ is a symmetric Newton polygon with $\np(G_0) \preceq \nu$,  then
there exists a principally quasi-polarized 
deformation of $G_0$ with generic Newton polygon
$\nu$.

\subsection{Monodromy of $p$-divisible groups}
\label{subsecmono}

As in the introduction, let $H_{r/s}$ be the $p$-divisible group with
Dieudonn\'e module $F^{s-r} = V^r$.  Dieudonn\'e and Manin have shown
that if $G$ is a $p$-divisible group over an algebraically closed
field $k$, then there exists an isogeny
\begin{equation}
\label{diagmanin}
\xymatrix{
\oplus_{\lambda\in\rat} H_{\lambda}^{\oplus m_\lambda} \ar[r] & G.
}
\end{equation}
Now suppose $G$ is a $p$-divisible group over an arbitrary field $K$ of
characteristic $p$.  (For ease of exposition below, we will always
assume that $K$ contains the algebraic closure $\bar\ff_p$ of the prime
field.)  In general, the slope filtration \eqref{eqslopefilt} only
splits, even up to isogeny, after passage to the perfect closure
$K^\perf$ of $K$ \cite{zinksf}.
Moreover, even if $K$ is perfect, an isogeny \eqref{diagmanin} need
not exist over $K$.  The field of
definition of such an isogeny is a measure of the complexity of the
$p$-divisible group.  Henceforth let $\gal(K) = \gal(\bar K/K^\perf)$;
it is canonically isomorphic to $\gal(K^\sep/K)$.

Let $\lambda= \lambda_i$ be one of the slopes of $G/K$, in the sense
that $m_\lambda > 0$.  Then  $\aut(H_{\lambda})$ acts on
$ V^{\lambda}= \hom(H_{\lambda}, (G_i/G_{i+1}))\otimes \rat$ on the right. 
Moreover, 
 $\gal(K)$ acts on $V^{\lambda}$ on the left; 
$\tau\in\gal(K)$
takes $f\in \hom_{\bar
  K}(H_{\lambda, \bar K}, (G_i/G_{i-1})_{\bar K})$ to $\tau\comp f \comp
\tau\inv$. The actions of $\gal(K)$ and $D_{\lambda}$ commute,
and we obtain a representation $\rho: \gal(K)
\rightarrow \aut_{D_{\lambda}}(V^{\lambda})$. 

\begin{definition}\label{defmono}
We call this the $\lambda$-monodromy of $G$, and the image of
$\gal(K)$ in $\aut((G_i/G_{i-1})_{\bar K})$ the $\lambda$-monodromy
group of $G$.  If $R$ is a complete local domain and $G/R$ is a
$p$-divisible group, the $\lambda$-monodromy of $G/R$ is that of
$G\cross{\Frac R}$.
\end{definition}

Given a choice of isomorphism $V^{\lambda}
\rightarrow D_{\lambda}^{m_{\lambda}}$ we obtain a representation
in $\gl_{m_{\lambda}}( D_{\lambda})$.
  Our goal is to  
show that the monodromy  differs little  
from $\gl_{m_{\lambda _i}}( \mathcal{O}_{\lambda})$. 
We say that two subgroups of $\gl_{m_\lambda}(D_{\lambda})$
are commensurable if  there is a single 
subgroup of finite index in both  groups. 

\begin{definition}\label{deflarge}  
We call a subgroup of
$\gl_{m_\lambda}(D_\lambda)$ {\em large} if it is commensurable with
$\gl_{m_\lambda}(\calo_{H_\lambda})$.
\end{definition}
Even though the $\lambda$-monodromy group depends on the choice of
isogeny \eqref{diagmanin}, of isomorphism $V^{\lambda}\rightarrow
D_\lambda^{m_\lambda}$, and of slope-$\lambda$ test object, we will
see below that having large $\lambda$-monodromy is independent of all
of these choices.

It will often be more convenient to calculate monodromy in the
category of $F$-lattices.  Assume $K$ is a perfect field of characteristic $p$ and let
$\sigma$ denote the Frobenius on $W(K)$. By  
an $F$-lattice we mean
a free, finitely generated $W(K)$-module with an injective
$\sigma$-linear operator $F$.  A Dieudonn\'e module over $K$ gives us
an $F$-lattice by forgetting the action of $V$.  We say that
$F$-lattices $M_1$, $M_2$ are isogenous if there is an $F$-equivariant
map $M_1 \ra M_2$ with $W(K)$-torsion kernel and cokernel.  If $M_1$, $M_2$ are two
 Dieudonn\'e modules over $K$ and if $\til M_i$ is an $F$-lattice which 
is isogenous to $M_i$ as
$F$-lattice for $i = 1, 2$, then\[
\hom_{\dieu}(M_1,M_2)\tensor\rat \iso \hom_F(\til M_1, \til
M_2)\tensor\rat,
\]
where the left-hand side denotes homomorphisms of Dieudonn\'e modules
and the right-hand side means homomorphisms of $F$-lattices.  Somewhat
more precisely, we have \cite[IV.1]{demazure}:

\begin{lemma} 
\label{lemftod}
Let $M_1$ and $M_2$ be Dieudonn\'e modules over a perfect field $K$
which are isogenous as $F$-lattices.  Then $M_1$ and $M_2$ are
isogenous as Dieudonn\'e modules, and $\hom_{\dieu}(M_1, M_2)$ has finite index 
in $\hom_F(M_1, M_2).$
\end{lemma}

Therefore, in order to compute $\lambda$-monodromy, we may work in the
category of $F$-lattices, rather than the category of Dieudonn\'e
modules. If $M$ and $N$ are two $F$-lattices, we say two  subgroups of
$\hom_F(M, N) \otimes \rat$ are commensurable if there is a single  
subgroup  of finite index in both. 

\begin{lemma}
\label{lemisogcommens}
For $i = 1,2$, let $M_i$ and $N_i$ be $F$-lattices over $K$.  If $M_1$
and $M_2$ are isogenous, and if $N_1$ and $N_2$ are isogenous, then
$\hom_F(M_1,N_1)$ and $\hom_F(M_2,N_2)$ are commensurable, as are
$\aut_F(M_1)$ and $\aut_F(M_2)$.
\end{lemma}
\begin{proof}
  If $M$ and $N$ are $F$-lattices, then $\hom(M,N)$ is naturally a
  summand of $\End(M\oplus N)$.  Therefore, for the first claim it
  suffices to prove that, for any pair of $F$-lattices $M_1$ and
  $M_2$, an isogeny $\phi:M_1 \ra M_2$ identifies an open subgroup of
  $\End(M_1)$ with an open subgroup of $\End(M_2)$.

Now, $\phi$ induces an isomorphism $M_1 \tensor\rat \iso
M_2\tensor\rat \iso V$; we view $\cale_i := \End(M_i)$ as a subgroup
of $\End(V)$.

Let $\cale_i(n) = p^n \cale_i = \st{\alpha\in \End(M_i) : \alpha(M_i)
  \subseteq p^n M_i}$.  Each $\cale_i(n)$ has finite index in
$\cale_i$.

Suppose that $p^n M_2 \subset M_1 \subset p^{-n}M_2$.  
Then $\alpha\in \cale_2(2n)$ maps
$M_1$ to itself.  Therefore, $\cale_2(2n)\subset \cale_1$.  In
particular, $\cale_2(2n) \subset \cale_1\cap \cale_2 \subset \cale_2$,
so that $\cale_1 \cap \cale_2$ has finite index in $\cale_2$.  After
taking an isogeny $M_2 \ra M_1$, we similarly see that $\cale_1 \cap \cale_2$
has finite index in $\cale_1$.

To prove the second claim we use $\phi:M_1 \rightarrow M_2$ to view
$\aut(M_1)$ and $\aut(M_2)$ as subgroups of $\aut(M_1 \otimes \rat)$.
Let $A=\{g \in \aut(M_1)\ | (g-1)M_1 \subseteq p^nM_1 \}$ 
for sufficiently large $n$. Then
$A \subseteq \aut(M_1) \cap \aut(M_2)$,  and $A$ has 
finite index in both $\aut(M_1)$ and $\aut(M_2)$.
\end{proof}

Below, for $\lambda = r/s$ we will need to consider the $F$-lattice
\begin{align}
\label{defnlambda}
N_{\lambda, F} &= W(\ff_p)[p^{1/s}][F]/(F-p^\lambda).
\end{align}
This $F$-lattice is isogenous to the
$F$-lattice obtained from the Dieudonn\'e module of $H_{\lambda}$.
Note that $N_{\lambda, F}$ admits an endomorphism $\varpi$ that maps
(the equivalence class of) $1$ to $p^{1/s}.$ For a field $K$
containing $\ff_q$, with $q=p^s$, the endomorphism ring of $N_{\lambda,
  F}(K)$ is $\calo_\lambda = W(\ff_q)[\varpi]$, an order in the
division algebra over $\rat_p$ with Brauer invariant $\lambda$.   Note
that $\varpi^s = p$, and that if $x\in W(\ff_q) \subset \calo_\lambda$
then $x\varpi = \varpi x^\tau$, where $\tau = \sigma^r \in
\aut(W(\ff_q)/\integ_p)$.  In fact, since $\gcd(r,s) = 1$, $\tau$ is a
generator of $\aut(W(\ff_q)/\integ_p)$.
The
automorphism group of $N_{\lambda, F}$ is
\begin{equation}
\label{eqdefcalg}
\calg_\lambda = W(\ff_q)[\varpi] \units;
\end{equation}
we consider the structure of $G_\lambda$ in Section
\ref{subsecsubgroup}.

\section{A special case}

In this section we prove a special, but crucial, case of our main result.  
Let $G_0$ be a local $p$-divisible group over an algebraically closed 
field $k$. We assume that the $a$-number of $G_0$ is one,  that
$\lambda=r/s$ is positive and less than any slope of $G_0$, that $\lambda
\ne (s-1)/s $ for any natural number $s$, and that $\lambda$ 
is attainable from $G_0$.
 Let $G^\univ/R^\univ$ be the  universal deformation of $G_0$ over the
 universal deformation ring $R^\univ$.
 Given a Newton polygon  $\nu$ with $\nu \succeq \np(G_0)$,
 by \cite{katzsf} there is a radical  ideal $J = J_\nu$ of $R^\univ$
characterized as follows: if   $R_1$ is  a complete local domain
and $G_{R_1}$ a deformation of $G_0$ whose Newton polygon
is on or above $\nu$, then the deformation $G_{R_1}$ is induced by a map
$R^\univ/J  \rightarrow R_1$; and if $x \in \spec(R^\univ/J)$, then
$\np(G_x) \preceq \nu$. Let $\nu=\np(*)$, the
Newton polygon obtained by adjoining $(s, r)$ to the 
Newton polygon of $G_0$. For this  choice of $\nu$ 
set $R = R^\univ/J$ and $G= G^\univ_R$.   (We given an explicit
description of $R$ in \eqref{diagdefdeform}.)
\begin{lemma} 
\label{lembasiccase}
Let $G_0/k$ be a local $p$-divisible group over an algebraically closed
field, with $a$-number $a(G_0) := \dim_k \hom(\aalpha_p, G_0) = 1$.
Suppose that $\lambda$ is a positive rational number strictly smaller
than any slope of $G_0$, that $\lambda\not = (s-1)/s$ for any natural number $s$, and that $\lambda$ is
attainable from $G_0$.  Let $G_R$ be the deformation described above.
Then the deformation $G_R$ of $G_0$ has large $\lambda$-monodromy.
\end{lemma}
(Note that by hypothesis $0 < \lambda < 1$ and $s \ge 3$.)
Let $K$ denote the fraction field of $R$. 
The hypotheses on $G_0$ and $\lambda$ force $a(G_K)=1$.  Therefore,
much information about $G_K$ is encoded in the (noncommutative)
characteristic polynomial of its Frobenius operator.  In Section
\ref{subsecdemazure}, we use a result of Demazure on factorization of
such polynomials to give a method for computing the lowest-slope
monodromy of a $p$-divisible group.  In Section \ref{subsecsubgroup}
we collect some remarks on the structure of $G_\lambda$, the ambient
group for the $\lambda$-monodromy group.

After reviewing how deformations of $G_0$ are described by displays
(Section \ref{subsecdisplay}), we analyze the  deformation
$G=G^\univ_R$ of $G_0$ in which $\lambda$ appears with multiplicity one
(Section \ref{subsecconstruct}).  We use the method of Demazure to explicitly
show that certain graded pieces of the $\lambda$-monodromy of $G$ are
maximal, and thus that the $\lambda$-monodromy of $G$ is large.

\subsection{A Lemma of Demazure}
\label{subsecdemazure}
In his book on $p$-divisible groups,
Demazure proves the following result \cite[Lemma IV.4.2]{demazure} about 
polynomials over the
noncommutative ring $W(K)[F]$, where $Fa = a^\sigma F$ for $a\in
W(K)$.

\begin{lemma}
\label{lemdemazure}
Given a polynomial in $W(K)[F]$,
\begin{equation}
\label{demazureinput}
\chi(F) = F^n + A_1 F^{n-1} + \cdots + A_n,
\end{equation}
let $\lambda = \frac rs = \min_i(\frac{\ord_p
A_i}{i})$ with $\gcd(r,s)=1$.  Over $W(\bar K)[p^{\oneover s}]$ there
exists a factorization
\begin{align*}
\chi(F) &= \chi_1(F) \chi_2(F) \\
&= \chi_1(F) \cdot (F-p^\lambda)u,\label{eqdemfact}
\end{align*}
where the element $ v = u\inv \in W( \bar K)[p^{\oneover s}]$ satisfies 
\begin{equation*}
 v^{\sigma^n} + a_1 v^{\sigma^{n-1}} + \cdots + a_n = 0,
\end{equation*}
with $a_i = A_i p^{-i\lambda}$.
\end{lemma}

In the situation where $\lambda$ is the smallest slope of $M$ and
occurs with multiplicity one, we can exploit such a factorization to
compute the monodromy of $M$ in slope $\lambda$ as the Galois group of
an equation such as \eqref{demazureinput}.  If $u\in W(\bar K)$, by
$K(u)$ we mean the extension of $K$ generated by all Witt components
of $u$.  Similarly, if $u\in W(\bar K)[p^{1/s}]$, then $u$ may be
written as $u = \sum_{i=0}^{s-1} u_i p^{i/s}$, and by $K(u)$ we mean
the extension of $K$ generated by all Witt components of each of the
$u_i$.

\begin{lemma}
\label{lemcalcgalois}
Let $M$ be a Dieudonn\'e module over a perfect field $K$ generated by a
single element $e$ that satisfies 
\begin{align*}
\chi(F)e = (F^n + A_1 F^{n-1} + \cdots + A_n)e = 0,
\end{align*}
where $n$ is the rank of $M$ as $W(K)$-module.  Suppose that the slope
$\lambda = \min_i(\frac{\ord_p A_i}{i}) = r/s$ occurs in $M$ with multiplicity one.
Set $a_i = p^{-\lambda i}A_i$.  Then the monodromy group of $M$ in
slope $\lambda$ is induced by the action of $\gal(K(v)/K)$
on $u=v^{-1}$, where $v\in
W(\bar K)[p^{1/s}]$ satisfies
\begin{equation}
\label{eqstar}
v^{\sigma^n} + a_1 v^{\sigma^{n-1}} + \cdots + a_n = 0.
\end{equation}
\end{lemma}

\begin{proof}
By Lemma \ref{lemftod}, it suffices to prove the result for
$F$-lattices, rather than for Dieudonn\'e modules.   As $F$-lattice,
$M$ is isogenous to the $F$-lattice $M_\chi = W(K)[F]/\chi(F)$.  Let
$N_\lambda = \dieu_*(H_\lambda)$ be the Dieudonn\'e module of slope
$\lambda$ defined in the introduction.  As $F$-lattice $N_{\lambda,
  F}$ is isogenous to $N_\lambda$, so we can replace the computation
of $\hom_\dieu(N_\lambda(\bar K), M(\bar K))$ with the calculation of
$\hom_F(N_{\lambda, F}(\bar K), M_\chi(\bar K))$.  Since $M_\chi$ is
defined over a perfect field it is isogenous to the direct sum of two
lattices $M_\chi \sim M_\chi' \oplus M_\chi''$, where the only slope
of  $M_\chi'$ is $\lambda$, with multiplicity one, and where $\lambda$
is not a slope of $M_\chi''$.
Thus we can replace our
calculation with the calculation of $\hom_F(N_{\lambda, F}(\bar K),
M'_\chi(\bar K))$.

Our next step is to use Demazure's lemma \ref{lemdemazure} to make
explicit the action of $\gal(K)$ on $M'_\chi(\bar K)$.  This
lemma gives us a map of $F$-lattices
\begin{equation*}
\xymatrix{
M_\chi(\bar K) = W(\bar K)[F]/\chi(F) \ar[r]^-\psi  &
M_\lambda(\bar K) = W(\bar K)[p^{1/s}][F]/(F-p^\lambda)u.
}
\end{equation*}
Denote the implicitly given generator of $M_\chi$ by $e_\chi$ and set
$e_{M_\lambda} = \psi(e_\chi)$.  For $x \in W(\bar K)$ and $g \in
\gal(K)$ we have $(xe_\chi)^g = x^g e_\chi$.  We claim that
$\ker\psi$ is Galois invariant, and thus $(x e_{M_\lambda})^g = x^g
e_{M_\lambda}$.  The map $\psi$ induces a diagram of quasi-morphisms
\begin{equation*}
\xymatrix{
M_\chi(\bar K) \ar[r] \ar[d] & M_\lambda(\bar K) \ar[d] \\
M'_\chi(\bar K) \oplus M''_\chi(\bar K)  \ar[r]^-{\til \psi} &
M_\chi(\bar K)
}
\end{equation*}
where the vertical maps are quasi-isogenies.  Since $\ker\til\psi$ is
Galois invariant, so is $\ker\psi$.

Let $e_{N_\lambda}$ denote the implicitly given generator of
$N_{\lambda, F}$.  Then
\begin{equation*}
\xymatrix{
N_\lambda \ar[r] & M_\lambda(\bar K) \\
e_{N_\lambda} \ar@{|->}[r] & u e_{M_\lambda}
}
\end{equation*}
is an isomorphism of $F$-lattices over $\bar K$.  The action of
$\gal(K)$ on $\hom_F(N_\lambda, M)$ is therefore the action of
$\gal(K)$ on $u$, so that $\gal(K(u)/K)$ induces an action
commensurable with the slope $\lambda$-monodromy of $M$.

\end{proof}

 Choosing an isomorphism
$M_{\lambda}(\bar K) \cong N_{\lambda , F}$ allows us to view 
the monodromy representation as a map from $\gal(K)$
to the group $\calg_{\lambda}$ defined in \eqref{eqdefcalg}. We next
make some preparatory remarks on the structure of $\calg_\lambda$.

\subsection{Subgroups of $\aut(N_{\lambda,F})$}
\label{subsecsubgroup}

We continue to let $\lambda = r/s$, where $r \ge 1$ and $\gcd(r,s) =
1$.  While the hypotheses of Lemma \ref{lembasiccase} imply that $s\ge
3$, the results of the present subsection are valid for $s = 2$, too.
Let $q = p^s$.
Recall \eqref{eqdefcalg}  that
$\calo_\lambda = \End(N_{\lambda, F})$ has a 
presentation as $W(\ff_q)[\varpi]$, where $\varpi^s= p$ and there is a
generator $\tau$ of $\aut(W(\ff_q)/\integ_p)$ such that, for all $x\in
W(\ff_q)\subset \calo_\lambda$, $x\varpi = \varpi x^\tau$.
Let $\calg = \calg_0 = \calo_\lambda\units$, and let $\calg_i = 1 + \varpi^i \calo_\lambda
\subset \calg$.   For any subgroup $\calh\subset \calg$, let $\calh_i
= \calh \cap \calg_i$.    There is always a natural inclusion
\begin{equation}
\label{diagdefgquot}
\xymatrix{
\calh_i /\calh_{i+1} \ar@{^{(}->}[r] & \calg_i / \calg_{i+1} \iso
{\begin{cases}
\ff_q\units & i= 0 \\
\ff_q^+ & i \ge 1
\end{cases}
}
.
}
\end{equation}

In Section \ref{subsecconstruct}, we will construct a $p$-divisible group whose
$\lambda$-monodromy $\calh$
satisfies $\calh_i/\calh_{i+1} = \calg_i /
\calg_{i+1}$ for $i = 0$, $1$ and $s$.  
We will show (Lemma \ref{lemthelemma}) that such a group
is in fact equal to $\calg$. 

The structure of an order in a division algebra over a local field is
efficiently documented in \cite{platonovrapinchuk}.
We make the isomorphisms in \eqref{diagdefgquot} explicit now.
An element $\alpha \in W(\ff_q)$ can be written as
\[
\alpha =   \sum_{j\ge 0} p^j \ang{\alpha_j},
\]
where $\alpha_j \in\ff_q$ and $\ang{\alpha_j} = (\alpha_j, 0, \cdots,
0)\in W(\ff_q)$. 
Now, the order $\calo_\lambda$ is a finite, free $W(\ff_q)$-module; any $\beta\in
\calo_\lambda$ has a unique expression
\[
\beta = \sum_{j=0}^{s-1} \varpi^j \beta_j,
\]
with $\beta_j\in W(\ff_q)$.  Taking the expansion of each $\beta_j$ as
above and relabeling, we have
\[
\beta = \sum_{j\ge 0} \varpi^j \ang{\beta_j}.
\]

\begin{lemma}
\label{lemsdivn}
Suppose that $\calh_s/\calh_{s+1} = \calg_s/\calg_{s+1}$.  Then
$\calh_{ns}/\calh_{ns+1} = \calg_{ns}/\calg_{ns+1} \iso \ff_q$ for all
natural numbers $n$. 
\end{lemma}
\begin{proof}
The proof is by induction on $n$.  Each class in
$\calg_{(n+1)s}/\calg_{(n+1)s+1}$ is represented by $1+\varpi^{(n+1)s}
\ang{\alpha}$ for some $\alpha\in \ff_q$.  By the
inductive hypothesis, there exists an element
$1+\varpi^{ns}\ang{\alpha} + \varpi^{ns+1}\til\beta \in \calh_{ns}$ for
some $\til\beta\in \calo_\lambda$.  Then
\begin{align*}
(1+\varpi^{ns}\ang{\alpha}+\varpi^{ns+1}\til\beta)^p &=
1 + p(\varpi^{ns}\ang{\alpha} + \varpi^{ns+1}\til\beta) +
p^2\varpi^{ns}\til\gamma
\\
&\equiv 1 + \varpi^{(n+1)s}\ang{\alpha}  \bmod \calg_{(n+1)s+1}
\end{align*}
for some $\til\gamma\in \calo_\lambda$, since $\varpi^s = p$.
\end{proof}

\begin{lemma}
\label{lemthelemma}
Suppose that $\calh_i/\calh_{i+1} = \calg_i/\calg_{i+1}$ for $i = 0$,
$1$ and $s$.  Then $\calh/\calh_n = \calg/\calg_n$ for all natural
numbers $n$.
\end{lemma}
\begin{proof}
It suffices to show that for all $n\in\integ_{\ge 0}$,
$\calh_n/\calh_{n+1} = \calg_n/\calg_{n+1}$.  Again, we prove this by
induction on $n$.  If $s|(n+1)$, then $\calh_{n+1}/\calh_{n+2} \iso
\calg_{n+1}/\calg_{n+2}$ by Lemma
\ref{lemsdivn}.  Otherwise, by induction we may assume that $\calh_i/\calh_{i+1} = \calg_i/\calg_{i+1}$ for $i =
1$ and $i = n$.  A direct computation \cite[Lemma
1.1.8]{platonovrapinchuk} shows that
$[1 - \varpi  \ang x, 1 - \varpi^n \ang y] \equiv 1 + \varpi^{n+1}
(x^{\tau^n} y - y^{\tau} x) \bmod \calg_{n+2}$.  Since $s\ndiv (n+1)$,
{\em every} element of $\ff_q$ is of the form $x^{\tau^n} y -
y^{\tau}x$; the result follows.
\end{proof}

\subsection{Displays}
\label{subsecdisplay}

Throughout this section we use the ideas, results, and terminology of
Oort's paper \cite{oortnpfg}.    If $R$ is a ring of positive
characteristic, and if $x\in R$, let $\ang x = (x, 0, 0, \cdots) \in W(R)$.


Let $k$ be a perfect field and $M_0$ a Dieudonn\'e module over $k$.
Denote the dimension of $M_0$ by $d$, the codimension by $c$, and the
height by $h=d+c$. By a display of $M_0$ we mean a choice of $W(k)$-basis
of $M_0$, $\st{e_i : 1 \le i \le h}$, along with relations defining
$M_0$:
\begin{align}
\label{eqdefdisp}
Fe_i & =  \sum a_{ji} e_j & 1 \le i \le d \\
e_i & = V(\sum a_{ji} e_j) & d+1 \le i\le h
\end{align}
We often summarize this data in the matrix
\begin{align}
\label{eqdefdispmat}
\begin{pmatrix}
A & B \\ C & D
\end{pmatrix}
\end{align}
where $A  = (a_{ij})_{ 1 \le i \le d, 1 \le j \le d}$,  $B =
(a_{ij})_{1 \le i \le d, d+1 \le j \le h}$, 
$C = (a_{ij})_{d+1 \le i \le h, 1 \le j \le  d}$,  and
$D = (a_{ij})_{d+1 \le i \le h, d+1 \le j \le h}$.

We define certain subsets of $\integ\cross \integ$: 
\begin{align*}
S &= \st{ (i,j) :  1 \le i \le d, d \le j \le h}\\
S^\univ  & = \st{ (i,j) : 1 \le i \le d, d+1 \le j \le h}.
\end{align*}
The universal equicharacteristic deformation of $M_0$ is defined over the ring
\begin{equation}
\label{eqdefruniv}
R^\univ := k\powser{t_{ij} : (i,j) \in S^\univ}.
\end{equation}
Let $T$ be the $d\cross c$ matrix with entries $T_{ij}
= \ang{t_{ij}}$.
Then the universal deformation of $M_0$ is displayed by
\begin{align}
\label{eqdefuniv}
\begin{pmatrix}
A + T C & B + TD \\
C & D
\end{pmatrix}.
\end{align}
(The theory of displays over rings which are not necessarily perfect
is documented in \cite{zinkdisp}.)

If the $a$-number of $M_0$ is one, then one can choose a basis for
$M_0$ so that the display \eqref{eqdefdisp} becomes particularly
simple.  Indeed, if the $a$-number of $M_0$ is one, so $\dim_k(M_0/(F,V)M_0)$ is
one,  there exists \cite[2.2]{oortnpfg}
a display so that the matrix $(a_{ij})$ has the form
\[
\left(
\begin{array}{l|l}
\begin{matrix}
0 & & \cdots & 0 & a_{1,d} \\
1 & 0 & \cdots & 0 & a_{2,d} \\
0 & 1 & \cdots & \\
\vdots & & & & \vdots \\
0 & \cdots &&1& a_{d,d}
\end{matrix}
&
\begin{matrix}
a_{1,d+1} & \cdots &&& a_{1,h} \\
a_{2,d+1} & \cdots \\
\vdots & \\ 
\\
\cdots &&&&a_{d,h}
\end{matrix}\\
\hline
\begin{matrix}
0 & \cdots &&& &1 \\
\vdots &&&&&0 \\
\vdots &&&&&\vdots \\
\\
0 & \cdots &&&& 0
\end{matrix}&
\begin{matrix}
0 & & \cdots &  & 0 \\
1 & 0 & \cdots &  & 0 \\
0 & 1 & \cdots & \\
\vdots & & & & \vdots \\
0 & \cdots &&1&0
\end{matrix}
\end{array}
\right)
\]
with  $a_{1,h}$  a unit in $W(k)$. We call such a display a normal
form for $M_0$.  In this
case, Oort shows there is an explicit polynomial $\chi_0(F)$ so that the 
generator $e$ of $M_0$ satisfies
\[
\chi_0(F)e = (F^h - \sum_{x=0}^{h-1} A_x F^{h-x})e = 0
\]
with $A_x \in W(k)$.  This polynomial depends on the choice of normal
display.  In spite of this ambiguity we call $\chi_0(F)$ the
characteristic polynomial of $M_0$.  Oort shows that the Newton
polygon of $M_0$ equals the Newton polygon of $\chi_0$.  We write out
the formula for the coefficient $A_x$ in terms of the entries $a_{ij}$
of the normal display, where $(i,j) \in S$.

Define a map
\begin{equation*}
\xymatrix{
\integ^2 \ar[r]^{f = (x,y)} & \integ^2 \\
(i,j) \ar@{|->}[r]
 & (x(i,j),y(i,j)) = (j+1-i, j-d).
}
\end{equation*}
With this notation, Oort's Cayley-Hamilton lemma \cite[2.6]{oortnpfg} says that
\[
A_x = \sum_{(i,j): (i,j)\in S, x(i,j) = x} p^{y(i,j)} a_{ij}^{\sigma^{h-y(i,j) - d}}.
\]
Observe that, since $\sigma$ is additive, this formula is additive in
$a_{ij}$.

 
The display for the
universal deformation $M^\univ$ is in normal form, and 
hence it is determined by
the entries in the  positions $(i,j) \in S$.  Let $\delta$ be the
translation map $\delta(i,j) = (i,j+1)$.
A display in normal form
for $M^\univ$ is given by the matrix $(a_{ij}^\univ)$, where
\begin{equation*}
a_{ij}^\univ =
\begin{cases}
a_{ij} + \ang{t_{\delta(i,j)}} & \delta(i,j) \in S^\univ\\
a_{ij} & 1 \le i \le d, j = h
\end{cases}
\end{equation*}


While the coordinates $t_{ij}$ on the deformation space arise
naturally from our choice of normal display, they obscure the Newton
stratification of that deformation space.  We introduce a new set of
coordinates $\til u_{x,y}$ on $R^\univ$ adapted to Newton polygon
calculations, as follows.

Let $\calp = f(\delta\inv(S^\univ))$; it is the set of lattice points in the
parallelogram with vertices $(1,0)$, $(d,0)$, $(c,c-1)$ and
$(h-1,c-1)$ (Figure \ref{pfig}).  For $(i,j)\in S^\univ$, let $\til
u_{x(\delta\inv(i,j)),y(\delta\inv(i,j))} = t_{i,j}$.  Then there is a
canonical isomorphism
\begin{equation}
\label{eqreparam}
R^\univ =
k\powser{\til u_{x,y}: (x,y)  \in \calp},
\end{equation}
 and the characteristic polynomial of $M^\univ$ is
\begin{equation*}
\chi_\univ(F) = \chi_0(F) + \sum_{(x,y) \in \calp} p^y \ang{\til
  u_{x,y}}^{\sigma^{h-d-y}}F^{h-x}.
\end{equation*}
Note that this gives a formula for computing the characteristic
polynomial of Frobenius of any deformation of $M_0$.  Indeed, let
$\phi: R^\univ \ra R$ be a map to a complete local ring, necessarily
of positive characteristic.  Write $r_{x,y} = \phi(\til u_{x,y})$;
then the characteristic polynomial of Frobenius of the deformation
$M/R$ of $M_0$ is
\begin{equation*}
\chi(F) = \chi_0(F) + \sum_{(x,y) \in \calp} p^y \ang{r_{x,y}}^{\sigma^{h-d-y}}F^{h-x}.
\end{equation*}

\subsection{Construction of a deformation}
\label{subsecconstruct}

We continue to work with a Dieudonn\'e module $M_0$ with $a(M_0) = 1$
and a positive slope $\lambda$ smaller than any slope of $M_0$ but
still attainable from $M_0$.  The possibilities for $\lambda$ are
completely determined by Theorem \ref{grothconj}.  Write $\lambda =
r/s$ with $\gcd(r,s) = 1$.  For $\lambda$ to be attainable, it is
necessary and sufficient that $r 
< c$, and that the slope of the line segment from $(s,r)$ to $(h,c)$
satisfies $\lambda \le (c-r)/(h-s) \le 1$.  In particular, we may and
do assume that $s>r$ and that $s \le r+d$.

We recapitulate the method of Oort for constructing deformations, and
then obtain details about  the characteristic polynomial 
of the resulting deformed module. As in the beginning of this section,
let $\np(*)$ denote the Newton 
polygon obtained by adjoining
the point $(s,r)$ to the Newton polygon of $M_0$; that is, $\np(*)$ is
the convex hull of $\np(M_0)\cup\st{(s,r)}$. 
Since the Newton polygon of a Dieudonn\'e module
with $a$-number one is the same as the Newton polygon of its 
characteristic polynomial, we can control the Newton polygon of a
deformation of $M_0$ by examining 
the $p$-adic ordinals of the coefficients of its
characteristic polynomial.

Define
\begin{align*}
\calp(*) &\defeq \st{(x,y) \in \calp : (x,y)\text{ lies on or above
  }\np(*)} \\
R & \defeq  k\powser{u_{xy} : (x,y) \in \calp(*)}.
\end{align*}
We define a deformation of $M_0/k$ to $M/R$ by specializing the universal
deformation:
\begin{equation}
\label{diagdefdeform}
\xymatrix{
R^\univ \ar[r]^\phi & R \\
\til u_{xy} \ar@{|->}[r]&
{\begin{cases}
u_{xy} & (x,y) \in \calp(*) \\
0 & (x,y) \not \in \calp(*)
\end{cases}}.
}
\end{equation}
By \cite[2.6]{oortnpfg}, the Newton polygon of $M$ is indeed $\np(*)$. We can say more.
Any deformation of $G_0$ to a complete local domain such that the  
Newton polygon of the deformed $p$-divisible group lies on or above $\np(*)$
arises from this deformation. This, together with the fact that $R$ is a 
domain, characterizes this deformation. 

Set $a_x(M_0)=p^{-\lambda x} A_x(M_0)$, and 
$a_x(M)=p^{-\lambda x} A_x(M)$. By Lemma \ref{lemdemazure}, 
the monodromy group of $M$ in slope $\lambda$
is the Galois group generated  by the Witt components of any $v$ which
satisfies 
\begin{equation}\label{eqmonolambda}
v^{\sigma^h} - \sum_{x=0}^{h-1}\left( a_x(M_0) + \sum_{y: (x,y) \in
  \calp(*)}p^{y-\lambda x} \ang{u_{x,y}}^{\sigma^{h-d-y}}\right) v^{\sigma^{h-x}} = 0.
\end{equation}

To prove that this Galois group is large, we need control over some
(in fact, three) of the terms which appear in \eqref{eqmonolambda}.
For any nonnegative integer $j$, let
\begin{align*}
\calp(j) &= \st{ (x,y) : (x,y) \in \calp(*), y - \lambda x = j/s} \\
\bar\calp(j) &= \cup_{i \le j} \calp(j).\\
\end{align*}
We will see in Section \ref{subseccalcgal} that partitioning
$\calp(*)$ as $\cup \calp(j)$ 
corresponds to keeping track of terms with $p$-adic valuation $j/s$.
In the present section, our goal is to show that $\calp(0)$,
$\calp(1)$ and $\calp(s)$ are nonempty.

\begin{lemma}\label{lemcalp0}
$\calp(0) = \st{(s,r)}$.
\end{lemma}
\begin{proof}
Certainly, $(s,r)$ is in $\calp(0)$.  Let $(t,u)$ be any lattice point
with $t \ge 0$ and $\lambda u = t$.  If $t < s$, then $\gcd(r,s)\not =
1$, contradicting our original hypothesis on the representation
$\lambda = r/s$.  If $t > s$ and $(t,u) \in \calp$, then the
sub-Dieudonn\'e module with slope at most $\lambda$ would have height
greater than $s$, which contradicts the definition of $\np(*)$.
\end{proof}

\begin{figure}
\setlength{\unitlength}{20pt}
\centerline{
\begin{picture}(10,5)
 \put(1,0){\makebox(0,0)[tl]{$(1,0)$}}
 \put(4,3){\makebox(0,0)[bl]{$(c,c-1)$}}
 \put(8,3){\makebox(0,0)[bl]{$(h-1,c-1)$}}
 \put(5,0){\makebox(0,0)[tl]{$(d,0)$}}
\put(1,0){\line(1,1){3}}
\put(4,3){\line(1,0){4}}
\put(1,0){\line(1,0){4}}
\put(5,0){\line(1,1){3}}
\put(0,0){\line(4,1){9}}
\put(4,1){\makebox(0,0)[b]{$(s,r)$}}
\multiput(0,0)(1,0){10}{\circle{0.01}}
\multiput(0,1)(1,0){10}{\circle{0.01}}
\multiput(0,2)(1,0){10}{\circle{0.01}}
\multiput(0,3)(1,0){10}{\circle{0.01}}
\multiput(0,4)(1,0){10}{\circle{0.01}}
\end{picture}
}
\caption{\label{pfig}{\em The parallelogram $\mathcal P$.}}
\end{figure}
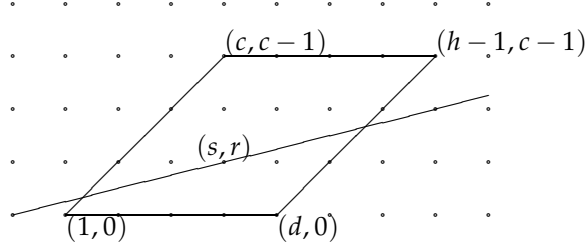

\begin{lemma} 
\label{lemcalp1}
If $(s,r)\not = (s,s-1)$, then $\calp(1)$ is nonempty.
\end{lemma}

\begin{proof}
Write $(s,r) = (a+b,a)$ with $\gcd(a,b)=1$, and consider Figure \ref{pfig}.  Since the segment
from $(a+b,a)$ to $(c+d,c)$ has slope at most one, it follows that $b
\le d$.  Since the Newton polygon is lower convex we may assume that
$a \le c-1.$

By hypothesis $b > 1$.   The line $y=(a/(a+b))x$ enters the
parallelogram $\mathcal P$ at $x=(a+b)/b$ and continues inside the
parallelogram at least
until the point $(a+b, a)$ since $d \ge b$.
For any integer $\alpha$
between $(a+b)/b$ and $a+b$, the point  with integer coordinates
directly on or  above
$(\alpha, \alpha \lambda)=
(\alpha, (a\alpha)/(a+b))$ is in $\mathcal P$. 
Thus it suffices to show
that there is an $\alpha$ in this range so that
\begin{align*}
 a \alpha/(a+b) &\equiv -1/(a+b) \mod{\integ} \\
\intertext{or}
 a \alpha &\equiv  -1 \mod{a+b},\\
\intertext{which is equivalent to }
b \alpha &\equiv 1 \mod{a+b}.
\end{align*}
There is a solution to this equation with $\alpha$ in the range
$1 \le \alpha < a+b$, since $\gcd(a,b)=1$. There is no solution
in the range $1  \le \alpha \le (a+b)/b$; hence there is a
solution in the range $(a+b)/b <  \alpha < a+b$.
\end{proof}

\begin{lemma}
\label{lemcalps}
If $(s,r)\not = (s,s-1)$,  then $\calp(s)$ is nonempty.
\end{lemma}

\begin{proof}
Since $r \le s-2$, $(s,r+1) \in \calp(*)$.
\end{proof}

\subsection{Calculation of Galois action}
\label{subseccalcgal}

Write $a_x=\sum p^{j/s}\ang{a_{x, j}}$ with $a_{x, j} \in k$.
By using Lemma \ref{lemcalp0} and recalling that $a_{x,0} = 0$,
we can
write our monodromy equation \eqref{eqmonolambda} as
\begin{equation}
\label{eqmonoax0}
v^{\sigma^h} - \ang{u_{s,r}}^{\sigma^{h-d-r}}v^{\sigma^{h-s}} -
\sum_{j \ge 1} \sum_x p^{j/s} v^{\sigma^{h-x}}
\left(\ang{a_{x,j}} + \sum_{y : (x,y) \in \calp(j)}\ang{u_{(x,y)}}^{\sigma^{h-d-y}}\right)
=0.
\end{equation}
(For given $x$ and $j$, if there is no $y$ such that $(x,y) \in
\calp(j)$ then we set the innermost sum equal to zero.)
Write $v = \sum p^{j/s} \ang{v_j}$, $K = \Frac R$, and $K_j = K(v_0,
\cdots, v_j)$. For a subset $P \subset \calp(*)$, let $k\powser{u_P} =
k\powser{u_{(x,y)}: (x,y) \in P}$, with field of fractions
$k\laurser{u_P}$.

\begin{lemma}
\label{lemcalcgal0}
$\gal(K_0/K) \iso \ff_q\units$.
\end{lemma}

\begin{proof}
Reduce equation \eqref{eqmonoax0} modulo $p^{1/s}$ and consider the
initial Witt component of $v$: it satisfies
\begin{equation}
\label{eqfirstwitt}
v_0^{p^h}-u_{r, s}^{p^{h-d-r}}v_0^{p^{h-s}}=0.
\end{equation}
Let $t$ be a nonzero solution to \eqref{eqfirstwitt}. It  satisfies
\begin{equation}
\label{eqnotwitt}
t^{p^h-p^{h-s}}=u_{r, s}^{p^{h-d-r}},
\end{equation}
and generates an extension of  separable degree $p^s-1$ 
over $k\laurser{u_{\calp(*)}}$  whose Galois group is $\ff_q\units$.
\end{proof}
We continue to let $t = v_0$.
The equality  \eqref{eqnotwitt} in the algebraic closure of $k\laurser{u_\calp}$ implies an equality of 
Witt vectors
\begin{equation}
\label{eqwittrel}
\frac{\ang{u_{s, r}}^{\sigma ^{h-d-r}}}{\ang{t}^{\sigma ^h-\sigma^{h-s}}}=1.
\end{equation}
In \eqref{eqmonoax0} set $v=\ang{t}w$ and divide by $\ang{t}^{\sigma ^h}$. Using \eqref{eqwittrel}, we obtain
\begin{equation*}
w^{\sigma ^h}-w^{\sigma ^{h-s}}-\ang{t}^{-\sigma ^h}
\sum_{j\ge 1} \sum_x p^{j/s} \ang{t}^{\sigma^{h-x}}w^{\sigma^{h-x}}
\left(
a_{x,j} + \sum_{y : (x,y) \in
  \calp(j)}\ang{u_{(x,y)}}^{\sigma^{h-d-y}}
\right)
= 0.
\end{equation*}

Write $w= \sum_{i=0} p^{i/s}\ang{w_i}$ with $w_0=1$.   Then our equation becomes
\begin{multline*}
(\sum_i p^{i/s}\ang{w_i})^{\sigma^h} - (\sum_i
p^{i/s}\ang{w_i})^{\sigma^{h-s}} \\
- \ang{t}^{-\sigma ^h}\sum_{j\ge 1} 
\sum_x
p^{j/s} 
\left(
\ang{a_{x,j}}+
\sum_{y : (x,y) \in \calp(j)}
\ang{u_{(x, y)}}^{\sigma^{h-d-y}}\right)
\left(\ang{t}^{\sigma^{h-x}} 
\sum_{i \ge 0} p^{i/s}
\ang{w_i}^{\sigma^{h-x}}\right) = 0
\end{multline*}
or, regrouping, 
\begin{multline*}
(\sum_\ell p^{\ell/s}\ang{w_\ell})^{\sigma^h} - 
(\sum_\ell p^{\ell/s}\ang{w_\ell})^{\sigma^{h-s}} \\
-\ang{t}^{-\sigma^h} \sum_{\ell\ge 1} p^{\ell/s} \sum_{j= 1}^\ell
\sum_x \left( \ang{a_{x,j}}
+ \sum_{y : (x,y) \in \calp(j)} \ang{u_{x,y}}^{\sigma^{h-d-y}}
\right)
\left(\ang{t}^{\sigma^{h-x}}\ang{w_{\ell-j}}^{\sigma^{h-x}}\right)=0.
\end{multline*}

Inductively, for $1 \le \ell \le s$ we have
\begin{equation}
\label{eqwell}
w_\ell^{p^h} - w_\ell^{p^{h-s}} -
t^{-p^h}\sum_{j=1}^\ell
\sum_x \left(
a_{x,j} + \sum_{y : (x,y) \in \calp(j)} u_{(x,y)}^{\sigma^{h-d-y}}
\right)
(t^{p^{h-x}}w_{\ell-j}^{p^{h-x}})=0.
\end{equation}
Since we have equality of Witt vectors in equation \eqref{eqwittrel}
the lifts $\ang{w_j}$ for $1 \le j < \ell$  solve the Witt 
equation modulo $p^{(1+s)/s}$. The  higher-order terms are irrelevant in
our calculations since $\ell \le s$.
Moreover, $w_\ell$ is defined in terms of $w_0, \cdots, w_{\ell-1}$ and
coordinates $u_{x,y}$ for $(x,y) \in \bar\calp(\ell)$, so that $w_\ell$ is
algebraic over $k\laurser{u_{\bar\calp(\ell)}}$.

Let $\calq(\ell) = \bar\calp(\ell) - \st{(s,r)}$.  With
this notation, $w_1, \cdots, w_{\ell-1}$ are algebraic over
$k\laurser{u_{\calq(\ell)}}\laurser t$.

Consider the extension 
$k\laurser{u_{\calq(\ell-1)}}\laurser{t} \subset 
k\laurser{u_{\calq(\ell-1)}}\laurser{t}(w_1, \cdots w_{\ell-1})$.
Since $k\laurser{u_{\calq(\ell-1)}}\powser{t}$ is a discrete valuation ring,
we can write this extension as
$$
k\laurser{u_{\calq(\ell-1)}}\laurser{t} \subset
L_{\ell-1}\laurser{t} \subset
L_{\ell-1}\laurser{t_{\ell-1}},
$$
where $L_{\ell-1}$ is an extension of the residue field 
$k\laurser{u_{\calq(\ell-1)}}$,  and the 
second extension is totally ramified with uniformizing parameter $t_{\ell-1}$.

The equation for $w_\ell$ is defined over
$$
k\laurser{u_{\bar{\calp}(\ell)}}(w_\ell, \cdots w_{\ell-1}) \subset
L_{\ell-1}\laurser{u_{\calp(\ell)}}\laurser{t_{\ell-1}},
$$
and has the form 
$$
w_\ell^{p^h}-w_\ell^{p^{h-s}}-
t^{-p^h}(\sum_{(x, y) \in \calp(\ell)} u_{x, y}^{p^{h-d-y}}t^{p^{h-x}})-B=0
$$ 
with $B \in L_{\ell-1}\laurser{t_{\ell-1}}.$
Write $A=t^{-p^h}(\sum_{(x, y) \in \calp(\ell)} u_{x,y}^{p^{h-d-y}}t^{p^{h-x}}).$

\begin{lemma}\label{lemnewvar} Assume $1 \le \ell \le s$. 
If there exists some $(x_0,y_0) \in \calp(\ell)$, then
the separable degree of the extension $L_\ell/L_{\ell-1}$ is greater than or equal 
to $p^s$. 
\end{lemma}

\begin{proof}

 Let $z=w^{p^{h-s}}$, then $z$ is a root of  $X^{p^s}-X=A+B$. 
This equation is separable; we show it is irreducible. By
Lemma \ref{lemappirred}, it suffices to show that it is impossible to write
$A+B=f_H(x)$, where $H$ is a nontrivial subgroup of $\ff_{p^s}$,
$f_H(x) = \prod_{a\in H}(x-a)$, and $x\in L_{\ell-1}\laurser{u_{\calp (\ell)}} \laurser{t_{\ell-1}}$. 

To show that $f_H(X)=A+B$ has no solution we use Lemma \ref{lemappnosol}.
Note that the field $k$ in the appendix corresponds to the field 
$L_{\ell-1}$ here; the variable $t$ in the appendix corresponds to $t_{\ell-1}$ here;
and the variables $z_1, \cdots z_e$ in the appendix correspond to 
$u_{x, y}$ with $(x, y) \in \calp(\ell)$.
\end{proof}

\begin{lemma}
\label{lemcalcgalpieces}
Suppose $r\not = s-1$.  Then $\gal(K_0/K) \iso \ff_q\units$, while
$[K_1:K_0] \ge p^s$ and $[K_s:K_{s-1}] \ge p^s$.
\end{lemma}
\begin{proof}
The first claim is Lemma \ref{lemcalcgal0}, while the rest follows
immediately from Lemmas \ref{lemcalp1}, \ref{lemcalps} 
and \ref{lemnewvar}.
\end{proof}

The main result of this section now follows easily.

\begin{proof}[Proof of Lemma \ref{lembasiccase}]
We show that the deformation defined by \eqref{diagdefdeform} has large monodromy in
slope $\lambda$.  The monodromy group admits a filtration by subgroups
of the form $1+\varpi^n \calo_\lambda$.  By comparing Lemma
\ref{lemcalcgalpieces} with the cardinality of the graded pieces of
$\calg_\lambda$ \eqref{diagdefgquot}, the
monodromy group is maximal for graded pieces $0$, $1$ and $s$.  By
Lemma \ref{lemthelemma}, we conclude that the $\lambda$-monodromy group is all
of $\aut(N_{\lambda,F})$.
\end{proof}

\section{Main result}
\label{sectmain}

The goal of this section is to complete the proof of

\begin{theorem}
\label{maintheorem}
Let $G_0/k$ be a $p$-divisible group.  Let $\lambda$ be a rational
number which is not a slope of $G_0$ but is attainable from $G_0$.
Assume that $\lambda\not = (s-1)/s$ for any natural number $s \ge 2$.
Then there 
exists a smooth
equicharacteristic deformation $G_R$ of $G_0/k$ 
to a domain $R$ so that $G_R$ attains $\lambda$, 
and the monodromy group of $G_K$ in slope $\lambda$ is large.
\end{theorem}

The theorem follows from Lemma \ref{lembasiccase} and the 
three lemmas below.
Lemma \ref{lemanumber} removes the hypothesis on $a(G_0)$ from Lemma
\ref{lembasiccase}.
Lemma \ref{lemsmallestslope} completes the proof of the main 
theorem for positive $\lambda $ by reducing to the case where $\lambda$ is
strictly less than any slope of $G_0$. Lemma \ref{lemslopezero}
proves the theorem 
 in case $\lambda$=0. 

\begin{lemma}
\label{lemanumber}
Let $G_0$ be a local $p$-divisible group over an algebraically closed field $k$. Let
$\lambda \in [0,1], \lambda \ne (s-1)/s $, 
be a positive rational number attainable from $G_0$ which is strictly less than any 
slope of $G_0$.
Then there exists a deformation of $G_0$ to a 
$p$-divisible group
$G$ which attains $\lambda$  with multiplicity one and that has large 
$\lambda$-monodromy. 
\end{lemma}

\begin{proof}

Suppose that $M_0$ has dimension $d$ and codimension $c$.  Choose a
display for $M_0$ as in \eqref{eqdefdisp}; then the universal
deformation of $M_0$, defined over $R^\univ$ and denoted $M^{\univ}$, 
is given in
\eqref{eqdefuniv}.

Oort shows \cite{oorttexelnp} that there exists a complete local 
domain $R$ and elements
 $\overline{r_{ij}} \in R$ so that, if we write $r$ for the matrix
 $(\ang{\bar{r_{ij}}})$,   then
the Dieudonn\'e module $M_1$ displayed by 
$$
\begin{pmatrix}
A + r C & B + r D \\
C & D
\end{pmatrix}
$$
has the same Newton polygon as $M_0$, but $a(M_1(\Frac R))=1$.

Let $S_{ij} = \ang{s_{ij}} \in W(R\powser{s_{ij} : (i,j) \in S^\univ})$.
Let $S=(S_{ij})$, and  let $M_2$ be the Dieudonn\'e module over $ R\powser{s_{ij}}$ 
with display
\begin{equation}
\label{eqdispm2}
\begin{pmatrix}
A + rC + S C & B + r D + S D \\
C & D
\end{pmatrix}.
\end{equation}
The Dieudonn\'e module $M_2$ is a deformation of $M_0$. Indeed, the
map
\begin{equation*} 
\xymatrix{ 
R^\univ \ar[r]^{\phi_1} & R\powser{s_{ij}} \\
t_{ij} \ar@{|->}[r] & r_{ij}+s_{ij}
}
\end{equation*}
exhibits $M_2$ as a deformation of $M_0$. 

Let $K$ be the fraction field of $R$ and $\bar K$ its algebraic 
closure. The 
Dieudonn\'e module $M_2(K\powser{s_{ij}})$
defined by \eqref{eqdispm2} is visibly the universal deformation
of $M_1(K)$. 

 Let $\lambda$ be a  slope attainable from $M_0$ which is strictly less than any 
slope of $M_0$.   Write $\lambda =
r/s$ where $r$ and $s$ are relatively prime integers.  
 Let $\np(*)$ 
denote the lower convex hull of the Newton polygon of $M_0$
with  $(s,r)$ adjoined.  By hypothesis, $\np(*)$ has the same
beginning and end points as the Newton polygon of $M_0$.

By a result of Katz \cite[2.3.1]{katzsf}, the locus of $\spf R^\univ$
over which the Newton polygon of the deformation lies on or above
$\np(*)$ is Zariski-closed.  Indeed, he  shows  there exists an ideal
$J\subset R^\univ$ so that, if $\phi:R^\univ \ra R'$ is any extension
of scalars to a domain $R'$ inducing a Dieudonn\'e module, 
then $\phi(J)R'$ defines 
(set-theoretically) the locus of points where the
Newton polygon $M^\univ(R')$ lies on or above $\np(*)$.

We apply this to our situation. Consider the sequence of  
homomorphisms
\begin{equation*}
\xymatrix{
R^\univ \ar[r] & R\powser{s_{ij}} \ar[r] & K\powser{s_{ij}} \ar[r] & \bar K \powser{s_{ij}}.
}
\end{equation*}
The ideals $J$, $JR\powser{s_{ij}}$,
$JK\powser{s_{ij}}$, and $J\overline{K}\powser{s_{ij}}$ define the
loci where $M^{\univ}$, $M_2$, $M_2(K)$, and $M_2(\overline{K})$,
respectively, specialize to a Dieudonn\'e module with Newton polygon
on or above $\np(*)$.

 Let $Y=R\powser{s_{ij}}/JR\powser{s_{ij}}$. Then $M_2(Y)$ is a 
deformation  of $M_0$ whose Newton polygon is equal to $\np(*).$
Our goal is  to show that there is a point 
of $\spf(Y)$ so that $M_2(Y)$ restricted 
to that point has generic Newton polygon equal to $\np(*)$
and has large monodromy.

Let $Y_{\overline{K}}$ denote $Y \otimes_R \overline{K}$.
The $a$-number of $M_1(\bar K)$ is one.  Since $M_2(\overline{K}\powser{s_{ij}})$ is the universal deformation 
of $M_1(\overline{K})$, 
and since $J\overline{K}\powser{s_{ij}}$ defines the locus of 
deformations of $M_1(\overline{K})$  
 with Newton polygon on or
above $\np(*)$, we can apply  Lemma \ref{lembasiccase}
to conclude that  the Dieudonn\'e module 
$M_2(Y_{\overline{K}})$ has large monodromy.

Let $I$ be the kernel of the map $Y \rightarrow Y_K$ and set
$Y'=
Y/I$.   We will show that $M_2(Y')$ has large monodromy.
To see this let $Q_{\bar K}$ be the field of fractions of
$Y_{\bar K}$, and let $Q'$ be the field of fractions of $Y'$.  
\begin{equation*}
\xymatrix{
&Q_{\bar K} = \Frac Y_{\bar K} \\
Y_{\bar K} \ar@{-}[ur] &Q' = \Frac Y'\subset \Frac Y_K \ar@{-}[u]\\
Y'\subset Y_K \ar@{-}[u] \ar@{-}[ur] &\\
}
\end{equation*}

Let $N_\lambda$ be the slope-$\lambda$ test object; it is defined over
$\ff_p$. We can replace the $p$-divisible group associated
to $M_2$ over the field $Q'$ by its local part. Hence we proceed  assuming that
$M_2$ is local.  Let $L$ be any algebraic extension of $Q'$ such that there
exists a nontrivial map
\begin{equation*}
\xymatrix{
N_\lambda\ar[r]^-\phi & M_2( Q'\ldot L).
}
\end{equation*}
We need to show that the action of $\gal(L/Q')$  is large.

Now, $\phi$ induces a map $\phi_{\bar K}$:
\begin{equation*}
\xymatrix{
  N_\lambda \ar[r]^-{\phi_{\bar K}} & M_2(Q_{\bar K}\ldot L).
}
\end{equation*}
By Lemma \ref{lembasiccase}, the action of 
$\gal(Q_{\bar K}\ldot L/Q_{\bar K})$ is
large.  The diagram of fields
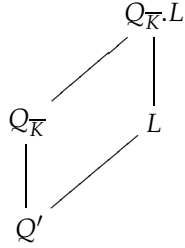
\begin{equation*}
\xymatrix{
&Q_{\bar K}\ldot L\\
Q_{\bar K} \ar@{-}[ur] &L \ar@{-}[u]\\\
Q' \ar@{-}[u] \ar@{-}[ur] &\\
}
\end{equation*}
yields an inclusion $\gal(Q_{\bar K}\ldot L/Q_{\bar K}) \inject
\gal(L/Q')$.  Therefore, $M_2(Y')$ has large $\lambda$-monodromy.
\end{proof}

\begin{lemma}
\label{lemsmallestslope}
Let $G_0$ be a $p$-divisible group over an algebraically closed field $k$.
Let $\lambda \in [0,1]$ be a positive rational number attainable from
$G_0$ which is not a slope of $G_0$. Then there exists a deformation of $G_0$ to a  complete local domain
with residue field $k$ that has large $\lambda$-monodromy.
\end{lemma}

\begin{proof}
 Since a $p$-divisible
group over a field always admits a slope filtration \eqref{eqslopefilt}, there exists a
filtration
\[
(0) = G_0^{(0)} \subseteq G_0^{(1)} \subseteq  G_0^{(2)} \subseteq
G_0^{(3)} = G_0
\]
such that:

\begin{itemize}
\item For $1 \le i \le 3$, the subquotient $H_0^{(i)} := G_0^{(i)}/G_0^{(i-1)}$ is
  a $p$-divisible group.
\item The slopes of $G_0^{(1)}$ are all less than $\lambda$; the
  slopes of $H_0^{(2)}$ are all greater than $\lambda$; the
slopes of $G_0^{(2)}$ are all less than $1$; and
  $H_0^{(3)}$ is \'etale.
\end{itemize} 
By definition of attainability $\lambda$ is attainable from $G_0$ if and only if it
is attainable from $G_0/G_0^{(1)}$.  Moreover, $\lambda$ is attainable
from $G_0/G_0^{(1)}$ if and only if it is attainable from
$H_0^{(2)}$; if not, a slope strictly larger than $1$ would appear in
a $p$-divisible group which attained $\lambda$.
 By Sublemma \ref{lemdeffilter}, a deformation $H^{(2)}/R$ of $H_0^{(2)}$
lifts to a deformation $G/R$ of $G_0$ as filtered $p$-divisible group.
Moreover, $\hom_{\bar K}(H_\lambda, G_{\bar K}) = \hom_{\bar
K}(H_\lambda, H^{(2)}_{\bar K})$.  Therefore, there exists a
deformation of $G_0$ with large monodromy in slope $\lambda$ if there
exists such a deformation of $H^{(2)}$.
\end{proof}

\begin{sublemma}
\label{lemdeffilter}
Let $(0) = G^{(0)}_0 \subseteq G^{(1)}_0 \subseteq \cdots\subseteq  G^{(r)}_0 = G_0$ be a
filtered $p$-divisible 
group over $k$ whose quotients $H^{(i)}_0 := G^{(i)}_0/G^{(i-1)}_0$ are $p$-divisible
groups such that $H^{(r)}_0$ is \'etale and $G^{(r-1)}_0$ is local.   Let $R$ be a local $k$-algebra, and for each $i$ let
$H^{(i)}$ be a deformation of $H^{(i)}_0$ to $R$.  Then there exists a
deformation $G^{(0)}\subseteq G^{(1)}  \subseteq \cdots \subseteq G^{(r)}$
of $G_0$ as filtered $p$-divisible group over $R$ such that
$G^{(i)}/G^{(i-1)} \iso H^{(i)}$ for each $i$.
\end{sublemma}

\begin{proof}
  This is essentially contained in \cite[2.4]{oorttexelnp}.  

First, we reduce to the case where the \'etale part of $G_0$ is
trivial.  Indeed, suppose $G^{(1)}_0\subset G^{(2)}_0 = G_0$ is an inclusion of
$p$-divisible groups so that $G^{(1)}_0$ is local and $G^{(2)}_0/G^{(1)}_0$ is \'etale.
Then $G_0$ admits a decomposition as a direct sum $G_0 \iso
(G_0/G^{(1)}_0) \oplus G_0^{(1)} = H_0^{(2)} \oplus H_0^{(1)}$. 
 If $H^{(i)}$ is a deformation of $H^{(i)}_0$ over
$R$ for $i=1,2$, then $H^{(1)} \subset H^{(1)} \oplus H^{(2)}$ is a 
suitable deformation of $G^{(1)}_0 \subset G^{(2)}_0$.

Henceforth, we assume that $G_0$ is local, so that $G^{(r-1)}_0 = G^{(r)}_0$.
Denote the
  dimension and codimension of $G^{(i)}_0$ by $d(i)$ and $c(i)$,
  respectively.  There is a $W(k)$-basis $x_1, \cdots, x_{d(r)}$,
  $y_1, \cdots, y_{c(r)}$ for $M = \dieu_*(G)$ so that $x_1, \cdots,
  x_{d(i)}$, $y_1, \cdots, y_{c(i)}$ is a $W(k)$-basis for
  $\dieu_*(G^{(i)}_0)$.  We describe the display
  $\begin{pmatrix}A&B\\C&D\end{pmatrix}$ with respect to this basis.

The matrix $A$ is block-upper-triangular, with blocks $A_1, \cdots
A_r$.  The matrix $A_i$ is square of size $d(i)-d(i-1)$; $B$, $C$ and
$D$ have an analogous structure.  Observe that $\begin{pmatrix} A_i &
  B_i \\ C_i  & D_i \end{pmatrix}$ is a display for $H^{(i)}_0$.  

Let $R$ be any local $k$-algebra, and suppose that a deformation
$H^{(i)}/R$ of $H^{(i)}/k$ is given for each $i$.  Such a deformation
is described by a display $\begin{pmatrix}A_i + T^{(i)} C_i & B_i +
  T^{(i)}D_i \\ C_i & D_i\end{pmatrix}$, where $T^{(i)}$ is a
$(d(i)-d(i-1)) \cross (c(i)-c(i-1))$ matrix with entries in $R$.  

Let $T$ be the block-diagonal matrix with blocks $T^{(1)}, \cdots,
T^{(r)}$, and use this to construct the display $\begin{pmatrix} A +
  TC & B+TD \\ C&D\end{pmatrix}$ over $R$.  This defines a
$p$-divisible group $G/R$.  Since each of $A+TC$, $B+TD$, $C$ and
$D$ is block-upper-triangular, the filtration of $G_0$ extends to a
filtration $G^{(0)} = 0 \subset G^{(1)} \subset \cdots \subset G^{(r)}
= G$ over $R$.  Moreover, the display of $H^{(i)}$ is identical
to that of $G^{(i)}/G^{(i-1)}$, so that these are isomorphic
$p$-divisible groups.
\end{proof}

We conclude this section by proving Theorem \ref{maintheorem} for
$\lambda = 0$.  Note that this gives a (somewhat anachronistic) proof
of Igusa's result \cite{igusa}.  We include the analogous result for
slope one, even though such a slope is not ``attainable'' in the
formulation of Definition \ref{defattainable}.

\begin{lemma}
\label{lemslopezero}
Let $R$ be a complete local ring with field of fractions $K$ and
residue field $k$.  Let $G/R$ be a $p$-divisible group, with special
and generic fibers $G_0$ and $G_K$, respectively.
\begin{alphabetize}
\item Suppose the multiplicities of slope zero in $G_K$ and $G_0$ are
$1$ and $0$, respectively.  Then the slope zero monodromy of $G_K$ has
finite index in $\integ_p\units$.
\item Suppose the multiplicities of slope one in $G_K$ and $G_0$ are
$1$ and $0$, respectively.  Then the slope one monodromy of $G_K$ has
finite index in $\integ_p\units$.
\end{alphabetize}
\end{lemma}

\begin{proof}
For (a), consider the representation
\begin{equation*}
\xymatrix{
\gal(K)  \ar[r] & \aut(\hom_{\bar K}(\mmu_{p^\infty}, G_{\bar K})) \iso
\integ_p\units.
}
\end{equation*}
Since the representation is continuous, the image of this map is
closed.  In fact, this image is actually infinite.  If not, then
over some finite extension $K'/K$ there would exist a nontrivial
homomorphism from $\mmu_{p^\infty,K'}$ to $G_{K'}$, which would
  necessarily extend \cite[1.2]{dejongtate} to the integral closure of
  $R$ in $K'$.  This contradicts the hypothesis that $0$ is not a
  slope of $G_0$.  In particular, the monodromy group contains some non-torsion element
$\alpha$.  Let $\beta=\alpha^{p-1}$.  Then there is some $n\ge 1$ so
that $\beta\equiv 1 \bmod p^n$ but $\beta\not \equiv 1 \bmod
p^{n+1}$.  Therefore, $\beta$ generates the group
$(1+p^n\integ_p)/(1+p^{n+1}\integ_p)$.  Since the Frattini subgroup of
$1+p^n\integ_p$ is $1+p^{n+1}\integ_p$, $\beta$ topologically
generates $1+p^n\integ_p$, a subgroup of finite index in
$\integ_p\units$. 

The proof of (b) is similar.  There is a subgroup $G^{(1)}_K \subset
G_K$ such that $G_K/G^{(1)}_K$ is geometrically isomorphic to
$\rat_p/\integ_p$, and the monodromy of $G$ in slope one is the same
as that of $G_K/G^{(1)}_K$.  We claim that this monodromy group is
infinite; as in part (a), this implies that the monodromy group has
finite index in $\integ_p\units$.

Indeed, suppose not.  Over some finite extension $K'/K$ there would
exist a nontrivial homomorphism $G_K \ra \rat_p/\integ_p$.  By
\cite[1.2]{dejongtate}, this would extend to a nontrivial homomorphism
$G_{R'} \ra \rat_p/\integ_p$, where $R'$ is the integral closure of
$R$ in $K'$.  This contradicts the hypothesis that $G_0$ is purely local.
\end{proof}

\section{Polarizations}

We now explain how to extend Theorem \ref{maintheorem} to the setting
of principally quasi-polarized $p$-divisible groups, since these are the ones which
arise in applications to abelian varieties.   For such a $p$-divisible
group the dimension and codimension are the same, and we set $c = d =
g$, $h = 2g$.

Recall that by a principally quasi-polarized (pqp) $p$-divisible group we mean a
$p$-divisible group $G$ equipped with a self-dual isomorphism $\Phi:G
\ra G^t$.  In this section, deformations will always be of pqp
$p$-divisible groups.  Recall that, given a deformation of a
$p$-divisible group, if a chosen quasi-polarization also deforms then it
does so uniquely.

Note that a slope $\lambda$ appears in the  Newton polygon of 
a pqp $p$-divisible group with same multiplicity as $1- \lambda$.
Moreover a pqp $p$-divisible group has large monodromy in 
slope $\lambda$ if and only if it has large monodromy in slope $1- \lambda$.

\begin{theorem}
\label{theoremmainpol}
Let $G_0$ be a principally quasi-polarized $p$-divisible group over an
algebraically closed field.  Suppose that $\lambda$ is not a slope of
$G_0$, that $\lambda\not = (s-1)/s$ for any natural number $s\ge 2$,
and that $\lambda$ is symmetrically attainable from $G_0$.
Then there exists a pqp deformation of $G_0$ which symmetrically
attains $\lambda$ with large $\lambda$-monodromy.
\end{theorem}

\begin{proof}
Since the proof is similar to that for $p$-divisible groups without
polarization, we limit ourselves to the points where the proofs
differ. Observe that slope $1/2$ is not symmetrically attainable. 
Since a slope $\lambda$ occurs in a pqp $p$-divisible group
if and only if $1-\lambda$ occurs, it suffices to  consider the 
case $\lambda < 1/2$.  The proof of Lemma \ref{lemslopezero} holds in
the presence of a principal quasi-polarization,  so we assume that $\lambda >0$.   

First, we show that it suffices to prove the theorem under the additional
hypothesis that $\lambda$ is strictly less than any slope of $G_0$.
We use techniques from Section 3 of Oort's article \cite{oorttexelnp},
in particular Lemmas 3.5, 3.7 and 3.8 there.  He proves that there
exists a filtration by $p$-divisible groups
\begin{align*}
(0) = G_0^{(0)} \subseteq G_0^{(1)} \subseteq  G_0^{(2)} \subseteq
G_0^{(3)} = G_0
\end{align*}
such that:
\begin{itemize}
\item For $1 \le i \le 3$, the subquotient $H_0^{(i)} := G_0^{(i)}/G_0^{(i-1)}$ is
  a $p$-divisible group.
\item The slopes of $G_0^{(1)}$ are all less than $\lambda$; the
  slopes of $H_0^{(3)}$ are all greater than $1-\lambda$;  the slopes
  of $G_0^{(2)}$ are all less than $1-\lambda$.
\item The filtration is symplectic, in the sense that the principal
  quasi-polarization $\Phi_0: G_0 \ra (G_0)^t$ induces isomorphisms
  $\Phi_0^{(i)}: H_0^{(i)} \ra (H_0^{(3-i+1)})^t$.
\end{itemize}
In particular, the pair $(H_0^{(2)}, \Phi_0^{(2)})$ is a pqp
$p$-divisible group.

We decompose the Newton polygon of $G_0$ into three parts.  The
initial, middle and final parts respectively arise from $H_0^{(1)}$,
$H_0^{(2)}$ and $H_0^{(3)}$.  Oort shows that a deformation of the pqp
$p$-divisible group $(H_0^{(2)}, \Phi^{(2)})$ to a ring $R$ lifts to a
deformation $(G,\Phi)$ of $(G_0, \Phi_0)$ so that the inclusion
$G_0^{(1)} \inject G_0$ extends to $G_0^{(1)} \cross R \inject G$.  In
particular, the Newton polygon of $G$ has the same initial and final
parts as $G_0$.  Thus, to prove the theorem it suffices to show that
there exists a deformation of $(H^{(2)}_0, \Phi^{(2)}_0)$ to a pqp
$p$-divisible group whose Newton polygon is the same as that of the
middle part of $G$ and which has large $\lambda$-monodromy.  In
particular, it suffices to prove the theorem under the hypothesis that
$\lambda$ is smaller than any slope of $G_0$.

Second, we show that we may assume that $a(G_0)= 1$.  By
\cite[Corollary 3.10]{oorttexelnp}, given a pqp $p$-divisible group
over a field, there exists a pqp deformation with the same Newton
polygon as the original group but with generic $a$-number one.  The
reduction argument of Lemma \ref{lemanumber} now applies.

We now assume that $a(G_0) = 1$, that $\lambda>0$, and that $\lambda$
is smaller than all slopes of $G_0$.  Since $G_0$ has no toric part
and is self-dual, it is local; we can use covariant Dieudonn\'e
theory.
Let $M_0 = \dieu_*(G_0)$, and let $g$ denote the dimension
and codimension of $M_0$.  By a result of Oort \cite[2.3]{oortnpfg} we can find a
$W(k)$-basis $\st{e_1, \cdots, e_{2g}}$ for $M_0$ so that the display
of $M_0$ is normal and the pairing takes the form
\begin{align*}
\ang{e_i,e_j}&=
\begin{cases}
1 & j = i+g \\
-1 & j  = i - g \\
0 & \abs{i-j}\not = g
\end{cases}.
\end{align*}
Define an involution on $S^\univ$
by $\Inv_M((i,j)) = (j-g,i+g)$.  The
universal pqp deformation $M^{\univ,\pol}$ of $M_0$ is defined over the ring
\begin{equation*}
R^{\univ,\pol} := R^\univ/(t_{ij}- t_{\Inv_M((i,j))}:
  (i,j)\in S^\univ).
\end{equation*}
We calculate the equation for $R^{\univ,\pol}$ with respect to the new
coordinates on $R^\univ$ introduced in \eqref{eqreparam}.  Let
$\Inv_{\np}(x,y) = (2g-x, g-x+y)$; we then have
\begin{equation*}
R^{\univ,\pol} = R^\univ/(\til u_{x,y} - \til u_{\Inv_{\np}(x,y)}).
\end{equation*}

Let $\np(*)_\pol$ denote the symmetric Newton polygon with endpoints
$(0,0)$ and $(2g,g)$ obtained from the Newton polygon of $M_0$ by
adjoining $(s,r)$.   It is the lower convex hull of $\np(M_0) \cup
\st{(s,r), \Inv_{\np}(s,r)}$. If $0 \le x \le s$, note that $(x,y) \in
\np(*)_\pol$ if and only if $\Inv_{\np}(x,y) \in \np(*)_\pol$.

The universal pqp deformation of $(M_0,\ang{\cdot,\cdot})$ with Newton
polygon on or above $\np(*)_\pol$ is defined over $R^\pol :=
R\tensor_{R^\univ} R^{\univ,\pol}$, where $R$ is the ring constructed
in Section \ref{subsecconstruct}; the associated $p$-divisible group
is the pullback of the tautological group $G^\univ/R^\univ$.

Furthermore, the analysis of section three goes through for
$G/R^\pol$, too.  Specifically, let $\calp(*)_\pol = \st{(x,y) \in
  \calp: (x,y)\text{ lies on or above }\np(*)_\pol}$, and let
$\calp(j)_\pol = \calp(j) \cap \calp(*)_\pol$.
Then Lemmas \ref{lemcalp0},
\ref{lemcalp1} and \ref{lemcalps} apply to $\calp(*)_\pol$, too, and
we have $\calp(0)_\pol = \st{(s,r)}$, and $\calp(1)_\pol$ and $\calp(s)_\pol$
are nonempty.  Therefore (cf. Section \ref{subseccalcgal})
$G^\pol/R^\pol$ has large $\lambda$-monodromy.
\end{proof}

\appendix

\section{Artin-Schreier equations}

We establish a criterion for  Artin-Schreier equations to be
irreducible. We use this to calculate our Galois group.

Let $q = p^s$, and let $G\subset (\ff_q,+)$ be a subgroup of order $p^N$.  Define
\[
f_G(X) = \prod_{\alpha\in G}(X-\alpha).
\]
Since $f_G(X+\beta)=f_G(X)$ for $\beta \in H$,
$f_G(X)=\sum_{j=0}N a_j X^{p^j}$. In particular,  $f_G$ is additive.

\begin{lemma}
\label{lemappirred} Let $K$ be any field containing an
$\ff_q$, and let $A\in K$.  Then
\[
F(X) = X^q - X -A
\]
is reducible if and only if $A = f_G(a)$ for some $a\in K$ and some
nontrivial subgroup $G\subseteq \ff_q^+$.
\end{lemma}
\begin{proof}
We write the polynomial as 
\begin{equation*}
F(X)=X^q-X-A=  f_{\ff_q}(X)-A.
\end{equation*}
Assume $F=\prod f_i$ is a non-trivial factorization of $F$
into irreducible monic factors. Let $y_1$ denote a root of $f_1$.
The roots of $F$ are $y_1+ \beta, \beta \in \ff_q$. Thus once we adjoin $y_1$
to $K$ we can split all the $f_i$ and thus the splitting fields of
the $f_i$ are all the same. There is a subgroup $H$ of $\ff_q$
so that $f_1(X) = \prod_{\alpha\in H} (X-(y_1+\alpha))$.
Since $f_H(X)-f_H(y_1)$
vanishes on the set $y_1+H$ and has leading coefficient $1$, we have
$f_1(X)=f_H(X)-f_H(y_1)$.
For $\beta \in \ff_q$, let $[\beta] = \beta + H$ denote the corresponding
equivalence class in $\ff_q/H$. Define
$$
f_{[\beta]}(X)=f_{\beta}(X)=f_1(X- \beta).
$$
This is independent of the choice of $\beta$ in $[\beta]$.  Since the
set of roots of $f_{\beta}$ is $y_1+\beta +H$, we
have $F=\prod_{[\beta] \in \ff_q/H} f_{[\beta]}$.
The constant term of $f_{\beta}=f_1(X-\beta)=f_H(X-\beta)-f_H(y_1)$ is
$$
f_H(-\beta)-f_H(y_1)=-(f_H(\beta)+f_H(y_1)).
$$
Thus
$$
A=\prod_{ [\beta] \in \ff_q/H}(-(f_H(\beta)+f_H(y_1))).
$$
The map $f_H$ is additive, maps $\ff_q$ to itself, and has kernel
exactly $H$. Let $L$ denote the image of $\ff_q$ under $f_H$. 
Then we can write
$$
A=f_L(-f_H(y_1)).
$$
Since $f_H(y_1)$ is the constant term in $f_1$, it is in $K$.
\end{proof}

\begin{lemma}
\label{lemappnosol}
Let $F(X)=X^{p^n}+c_{n-1}X^{p^{n-1}}+ \cdots +c_0X$ be an 
additive polynomial  with coefficients in a finite field.  
Let $k$ be any field containing the coefficients of $F$, and let $K =
k\laurser{z_1, \cdots z_e} \laurser t$.  Suppose 
\[
A = z_1^M( d_{-N}t^{-N} + d_{-N+1} t^{-N+1} + \cdots)+B  \in K
\]
where $M,N\in \nat$, $d_i \in k$, $d_{-N}\not = 0$, and $B \in k\laurser t$.
  Then $F(X) = A+B$ has no solution
in $K$.
\end{lemma}

\begin{proof} 
Possibly after enlarging $k$, we may and do assume that $d_{-N} = 1$.
We define an endomorphism $\pi$ of $K$ considered as $k$-vector
space.
Let $j \in \integ$, and $\alpha \in \integ^e$.
  Let
\begin{equation*}
\pi_j(z^{\alpha}) =
\begin{cases}
z^{\alpha} & \text{if } \alpha=(i, 0 \cdots 0) \text{and } \frac ij = \frac{M}{-N} \\
0 & \text{otherwise}\\
\end{cases},
\end{equation*}
Now extend $\pi_j$ to Laurent series by setting 
$\pi_j(\sum_{\alpha}a_{\alpha}z^{\alpha})=
\sum_{\alpha}a_{\alpha}\pi_j (z^{\alpha})$. 
For $x=\sum_ j A_jt^j$  $A_j \in k\laurser{z_1, \cdots, z_e}$, define 
$\pi(x)=\sum \pi_j(A_j)t^j$.

  Then $\pi$ is an
idempotent operator; $\pi\comp\pi = \pi$.  Moreover, $\pi$ commutes
with the $p^{th}$-power map; for $\alpha \in K$, $\pi(\alpha^p) =
(\pi(\alpha))^p$.  In particular, $\pi$ commutes with $F$ as $k$-linear endomorphisms of $K$.

This implies that if $F(X)=A+B$ has a solution in $x \in K$,  then 
$\pi(x)$ is a solution to $F(X)=\pi(A+B)=\pi(A)$.    Indeed, if $F(x)
= A+B$, then $F(\pi(x)) = \pi(F(x)) = \pi(A+B)$.  Since $\pi(A+B) =
\pi(A) = z^Mt^{-N}$, we are reduced to showing that $F(X)=z^Mt^{-N}$ has no 
solution of the form $\pi(x)$ for $x \in K$.

Write $e=\gcd(M,N)$, $m = M/e$, $n = N/e$, and $w =z_1 ^mt^{-n}$. The
$k$-linear space $\pi(K)$ can be identified with $k\laurser{w\inv}$.
We are reduced to showing that
$$
F(X)=w^e
$$
 has no solution in $k\laurser{w\inv}$. Since $F$ maps 
each of $k[w]$ and $k\powser{w^{-1}}$ to itself, we are further reduced to showing that
$F(X)=w^e$ has no solution in $k[w]$. But $F(x)$ never produces a monomial for
any $x \in k[w]$ of positive degree. 
Thus the equation $F(x)=w^e$ has no solution.
\end{proof}

\bibliographystyle{plain}
\bibliography{jda}

\end{document}